\documentclass[12pt]{amsart}   
\usepackage{amssymb,amscd,latexsym}
\usepackage{amsmath}

\newcommand{\transpose}{\hskip 2pt{}^{{}^{\rm T}}\hskip -3pt}

\newcommand{\rar}{\rightarrow}

\newtheorem{Theorem}{Theorem}
\newtheorem{Lemma}[Theorem]{Lemma}
\newtheorem{Corollary}[Theorem]{Corollary}
\newtheorem{Proposition}[Theorem]{Proposition}
\newtheorem{Definition}[Theorem]{Definition}
\theoremstyle{remark}

\newtheorem{Remark}[Theorem]{Remark}

\newtheorem{Example}[Theorem]{Example}

\def\qed{\hspace*{\fill} $\square$}
\def\QED{\qed}
\def\J{{\mathcal J}}
\def\O{{\mathcal O}}
\def\cV{{\mathcal V}}
\def\h{{\mathfrak h}}
\def\g{{\mathfrak g}}
\def\p{{\mathfrak p}}

\def\der#1{{\rm Der}(#1)}
\def\hder#1#2{{\rm Der}_{#1}(#2)}

\def\cc{\mathbb C}

\def\rank{{\rm rank}\,}
\def\grade{{\rm grade}\,}
\def\Ass{{\rm Ass}\,}

\def\spec#1{{\rm Spec}\, (#1)}

\begin{document}           

\title{Conormal modules via primitive ideals}
\author{Guangfeng JIANG}

\thanks{Supported by JSPS, LNSTC, NNSFC}

\address{{\it Permanent}:  {\rm Department of Mathematics,
 Jinzhou Normal University,
 Jinzhou, Liaoning, 121000, P. R. China}}
\email{jzjgf@mail.jzptt.ln.cn}
\address{{\it Till March 31, 2000}: {\rm Department of Mathematics,
 Faculty of Science,
Tokyo Metropolitan University, Minami-ohsawa 1-1, Hachioji-shi, Tokyo 193-0397,
Japan}}
\email{jiang@comp.metro-u.ac.jp}
\keywords{conormal module,  primitive ideal,
symbolic power, torsion}
\subjclass{Primary 13C12, 14B07; Secondary  32S05, 32S30}

\begin{abstract}
The main object of this note is to study the conormal module $M$ and the
computation of the  
second symbolic power $\bar{\g}^{(2)}$ of an ideal $\bar{\g}$ in the
 residue  ring $\O/\h$ of a polynomial
 ring $\O$ over a field of characteristic zero. 
The torsion part $T(M)$ of $M$ and  the torsion free module $M/T(M)$ 
 are expressed by the  primitive ideal of $\g$ relative to $\h$.
Two characterizations for  $M/T(M)$ to be free are proved.
Some immediate applications are worked out.

\end{abstract}
\maketitle

\section{Introduction}\label{introduction}

Let $\h\subset \g$ be two ideals of a polynomial ring $\O$ over a
field $k$ of characteristic zero. Let
$\bar{\g}:= \g/\h$ be the image of $\g$ in $\O/\h$ under the canonical
projection. The main interest of
 this article is to study the conormal module $M:=\bar{\g}/\bar{\g}^2$
 and the  second symbolic power $\bar{\g} ^{(2)}$ of $\bar{\g}$.
 The connection between $M$ and $\bar{\g} ^{(2)}$ is established by
  the primitive ideal of $\g$ (relative to $\h$), which was  introduced by
 Siersma-Pellikaan \cite{P1, P2} and generalized to relative version 
in \cite{J, JS}.

In general,  $\O/\h$ is not regular, so
 the $\O/\g$  module  $M$ is neither  free nor 
torsion free even if both $\h$ and $\g$ are complete intersections,
and the projective dimension of $\bar{\g}$
is not finite.
 Especially, no generating set of $\bar{\g}$ 
forms   a regular sequence. 
Then the following questions would be interesting.
\begin{itemize}
\item[a)]{ Find descriptions of 
 the torsion part $T(M)$ of $M$, calculate the length (when it is finite) of
$T(M)$;}
\item[b)]{ Find descriptions of   the torsion free module $N:=M/T(M)$ and
conditions on the   freeness of $N$.}
\end{itemize}
In commutative algebra, there is a question by Vasconcelos
 on how to compute effectively the
symbolic powers of an ideal in a residue ring of a polynomial ring 
\cite{Aron1, vas2}. It follows from  \cite{JS} that the second symbolic power
of $\bar{\g}$ is the image of the primitive ideal of $\g$  in the
residue ring. We give a precise expression of $\bar{\g}^{(2)}$ under
some assumptions.

What brought our attention to these questions is the studying of
functions with non-isolated singularities on singular spaces.
In general the ideals $\h\subset \g$ define two subvarieties
$X\supset \Sigma$ of $\cc ^n$ if $k=\cc$. The primitive ideal  
of $\g$ collects all the functions whose zero level hypersurfaces
pass through $\Sigma$ and are tangent to the regular part  $X_{\rm reg}$
of $X$ along $\Sigma \cap X_{\rm reg}$. If we  supply $X$ with the so
called {\it logarithmic stratification }\cite{Saito}, then the 
primitive ideal of $\g$ 
consists of exactly
all the functions from $\g$ whose stratified  critical loci on $X$ 
contain $\Sigma$ (cf. \cite{J}). Hence, locally
 the primitive ideal plays a  similar role  to the second power of the maximal
ideal of the local ring $\O_{\cc ^n, 0}$ in singularity theory.
In order to study the topology of the Milnor fibre $F_f$ of a function
$f$ with singular locus $\Sigma$, we use a good deformation 
(the Morsification) $f_s$ of $f$. This $f_s$ has relatively
simpler singularities than $f$.  The existence
of the good deformation and related invariants (both topological 
and algebraic) have close relationship with $M$, $T(M)$ and 
$N$. Roughly speaking, the freeness of $N$  implies the existence
of the good deformation \cite{J}. The length of  torsion module $T(M)$
(when it is finite) gives some
 information on how $\Sigma$ sits in $X$ (cf. \cite{JS1}).

Under some conditions, we answer the questions a) and b). More precisely,
after some descriptions of $T(M)$ and $N$,
we mainly prove the following (see also Remark~\ref{critera}) 

\vskip 1mm

\noindent{\bf Main  Theorem }\label{conclusion} {\sl
Let $\h\subset \g$ be  complete intersection ideals of a
 polynomial ring $\O$ over a field $k$ of characteristic zero. Let 
  $\spec{\O/\g}$ be reduced and connected,
and the Jacobian ideal  $\J(\h)$ of $\h$ be not contained in  any
minimal  prime of
 $\O/\g$.  The $\O/\g$-module  $N$ is free if and only if
there exists an $\O$-regular sequence $g_1, \ldots , g_n$
 generating $\g$, such that 
$$\int _{\h}{\g}=(g_1, \ldots , g_p)+(g_{p+1},  \ldots , g_n)^2,
$$
where $p:=\grade\h,\,  n:=\grade \g$.
}

\vskip 1mm

As an application, in  the last  section we study 
 lines on a variety with isolated complete intersection singularity. More
applications to the general deformation theory of non-isolated
singularities on singular spaces  will be given in the sequel papers.

\section{Primitive ideals}\label{primitive ideals}

Let ${\O}$ be a commutative ring and $k\subset {\O}$ be a subring.
Let  $\der {\O}$ denote the 
$\O$-module of all  the $k$-derivations of $\O$.
For an ideal $\h\subset \O$, define
$$\hder {\h} {\O}:=\{\xi \in \der {\O}\mid \xi (\h) \subset \h\}.$$

\begin{Definition}\label{primitive}\rm
Let $\h\subset \g\subset \O$ be ideals. The
{\it  primitive ideal of $\g$} relative to $\h$ is 
$$\int _{\h}\g:=
\{f\in \g\mid \xi (f)\in \g \text{ for any } \xi \in \hder {\h}{\O}\}.$$
\end{Definition}

\begin{Remark}\label{geometry}\rm
 This definition  is a generalization of \cite{P1, P2}, and was generalized
to higher order relative version in \cite{JS}.
 And under general assumptions (cf. Proposition~\ref{torsionpart}), 
it was proved in \cite{JS} 
that the  primitive ideal of $\g$ relative to 
$\h$ is the inverse image in $\O$ of the second symbolic power of
the quotient ideal $\g/\h \subset \O / \h$.
\end{Remark}

Some basic facts about the primitive ideals are collected 
in the following lemma.

\begin{Lemma}\label{property}\sl
 Let $\h\subset \g$ be ideals of $\O$,  a commutative noetherian 
$k$-algebra.
Then
\begin{itemize}
\item[(1)]{  $\int _{\h}\g$ is contained in $\g$ and 
 contains $\h$ and ${\g} ^2$; }

\item[(2)] {for any $\g _i\supset \h$ {\rm (}$i=1,2${\rm )},
$\int _{\h}\g _1 \cap \g _2=\int _{\h}\g _1 \cap \int _{\h}\g _1$;}

\item[(3)] {If $\O$ is a polynomial ring over a field of characteristic
 zero, $\left(\int _{\h}\g\right)/{\h}=\int _0(\g/\h)$, where $0=\h/\h$;}

\item[(4)]{If $\O$ is a $k$-algebra of finite type over a commutative
 ring $k$, and   $\h$ and  $\g$ have no embedded primes, then
 for any multiplicative set ${\mathcal S}\subset \O$, 
${{\mathcal S}}^{-1}\left(\int _{\h}\g\right)\simeq \int _{{{\mathcal S}}^{-1}\h}{{\mathcal S}}^{-1}\g$.}
\end{itemize}
\end{Lemma}

\begin{proof}
The first three statements follows immediately from the definition. 
To prove (4),
one may use the fact that ${{\mathcal S}}^{-1}\hder {\h}{\O} \simeq 
\hder {{{\mathcal S}}^{-1}\h}{{\mathcal S}^{-1}\O}$. See \cite{JS} for 
details.
\end{proof}

\section{Conormal module: the torsion part}\label{torsionsection}

Let $M$ be a module over a ring $\O$, and $S$ be the set of all the non-zero
 divisors of $\O$. Denote by $Q$ the total quotient ring of $\O$ with 
denominator set $S$, and by $M_{S}:=M\otimes Q$ the quotient module of $M$.
The kernel  of the canonical map $M\rar M_S$ is  denoted by $T(M)$
and  is called the torsion (part) of $M$. 

Let $\h \subset \g $ be ideals of the commutative noetherian ring $\O$.
Denote by $\bar{\g}:={\g}/{\h}$,  the quotient ideal in $\O/\h$. The
  $\O/\g$-module 
$M:={\bar{\g}}/{{\bar{\g}}^2} \simeq {\g}/{{\g}^2 +\h}$
is called the conormal  module of $\bar{\g}$.

\begin{Proposition}\label{torsionpart}\sl
Let $\O$ be a polynomial ring over a field $k$ of characteristic 0. Let 
$\h\subset \g $ be unmixed ideals of $\O$ with $\g$ radical,
such that the Jacobian ideal
 $\J(\h)$ of $\h$ is not contained in any minimal prime of $\O/\g$, then
$$T(M)=T:=\frac{\int _{\h}\g}{\g ^2 +\h}
\simeq \frac{\bar{\g}^{(2)}}{\bar{\g}^2}. $$
Consequently,  we have the following exact sequence
$$
0\longrightarrow T(M){\buildrel{\imath}\over{\longrightarrow}} M
{\buildrel{\pi}\over{\longrightarrow}} N\longrightarrow 0 
\eqno{(\ref{torsionsection}.1)}$$
where  $N:=\g / \int _{\h} {\g}$.
\end{Proposition}

\begin{proof}
In the following we use $\bar{\p}$ to denote the image of an ideal
 $\p$ of $\O$ in the residue  ring $\O/\h$.
We first prove that $T\subset T(M)$. Note that for any 
${\p}\in \cV(\g)\setminus ({\cV}(\J(\h))\cup {\cV}(\J(\g)))$, 
$T_{\p}=
\left(\int _{0}\bar{\g}_{\bar{\p}}\right)/{\bar{\g}^2_{_{\bar{\p}}}}$
is a  module over local ring $R:=\left(\O/\h\right)_{\bar{\p}}$.
Since $R$ and $R/\bar{\g}_{\bar{\p}}$ are regular,
 $\bar{\g}_{\bar{\p}}$ is generated just by a part of the regular system 
 of parameters $\{g_1, \ldots, g_d\}$.
(cf. \cite[Theorem 36]{matsumura}). We may also assume that
 $\g$ is prime. It follows from \cite{hochster}
that $\bar{\g}_{\bar{\p}}^2=\bar{\g}_{\bar{\p}}^{(2)}$, 
the second symbolic power of $\bar{\g}_{\bar{\p}}$.
By \cite{Seibt}, 
$\int _{0} {\bar{\g}_{\bar{\p}}}=\bar{\g}_{\bar{\p}}^{(2)}$.
Hence $T$ is annihilated by a power of $\J(\g)\J(\h)$.

For any $\bar{a}\in T(M)$, let $\bar{\beta}\in \O/\g$ be a
 non-zero divisor such that $\bar{\beta}\bar{a}=\bar{0}.$
By taking representative, we have $\beta a \in \g ^2+\h$. For any 
$\xi \in \hder {\h}{\O}$, we have
 $\xi (\beta a)\equiv \beta \xi (a) \mod  \g.$ Since 
$\xi (\beta a)\in \g$, hence $\xi (a) \in \g$, and $a \in \int _{\h}\g$.

By \cite{JS}, $\left(\int _{\h}\g\right)/\h=\bar{\g}^{(2)},$ hence
$T\simeq \bar{\g}^{(2)}/\bar{\g}^2$.
\end{proof}

\begin{Lemma}\label{Ass2}\sl
Let $\g$ be a radical ideal of a commutative noetherian ring
 $\O$. Suppose that $\g$ is generated by an  $\O$-regular sequence
$g_1, \ldots , g_n$. If $\g ^*:= (g_1, \ldots , g_t)+ (g_{t+1}, \ldots , g_n)^2$
for some $1\le t\le n$, then $\Ass (\O/\g ^*)=\Ass (\O/\g).$
\end{Lemma}
\begin{proof} We use the  fact: For two ideals  $J\subset I$ with
$I$ radical and $\sqrt[]{J}=I$, then
$\Ass (\O/J)=\Ass (I/J)\cup \Ass (\O/I).$

Since $\sqrt[]{\g ^*}= \g$,  we need to prove that 
$\Ass (\g/\g ^*)\subset \Ass (\O/\g). $
Let $\g=\bigcap\limits_{i=1}^r P_i$ be the unique minimal prime
decomposition of $\g$. Let $Q\in\Ass (\g/\g ^*) $. By definition,
there exist $0\ne \bar{x}_i\in \O/\g$ and $0\ne \bar{y}\in\g/\g ^* $
such that $P_i={\rm Ann}(\bar{x}_i)$ ($i=1, \ldots, r$) 
and $Q={\rm Ann}(\bar{y})$. Suppose that 
$q\in Q\setminus\bigcup\limits_{i=1}^rP_i\ne \emptyset$. Then $qy\in
\g ^*$, and $qx_i \notin \g$ ($i=1, \ldots, r$). Write
$qy\equiv
\beta _1 g_1+\cdots +\beta _{t}g_t  \mod  \g^2
\quad {\rm and } \quad y= \sum_{k=1}^n y_kg_k.$
Then
$qy- (\beta _1 g_1+\cdots +\beta _{t}g_t)
=(qy_1-\beta _1)g_1+\cdots +(qy_t-\beta
_t)g_t+qy_{t+1}g_{t+1}+\cdots + qy_ng_n
\equiv 0 \mod  \g ^2.$
It follows that $qy_{t+j}\in \g=\bigcap\limits_{i=1}^rP_i$.
 Since $q\notin P_i$ for all 
$i=1, \ldots, r$, $y_{t+j}\in \bigcap\limits_{i=1}^rP_i=\g$, which
implies that  $y\in \g ^*$, a contradiction. 
We may assume that $Q\subset P_1$.

On the other hand, since $\Ass (\g/\g ^*)\subset \Ass (\O/\g ^*)$
and ${\rm rad}(\g ^*)=\g$, we have $Q\supset \g=\bigcap\limits_{i=1}^rP_i. $
It follows from \cite[(1.11)]{AM} that $Q\supset P_i$ for some
 $1\le i\le r$, which implies that $i=1$ since
$\g=\bigcap\limits_{i=1}^rP_i$
is the minimal prime decomposition of $\g$.
\end{proof}

\begin{Proposition}\label{basic}\sl
Let $\h\subset \g$ be  ideals of $\O$, a polynomial ring over a
 field $k$ of characteristic zero. Assume that $\h$ is unmixed and 
 $\g$ is  a radical complete intersection ideal.
 Suppose that the Jacobian ideal
 $\J(\h)$ of $\h$ is not contained in any minimal prime of $\O/\g$.
If there exists a minimal generating set $\{g_1, \ldots, g_n\}$ of $\g$
 such that
\begin{itemize}
\item[(1)]{The images of $g_1, \ldots , g_t$ are contained in $T(M)$
  for some integer $t: \, 0\le t\le n$;}
\item[(2)]{For any prime  $\p\in \Ass (\O/\g)$, 
 the images of $g_{t+1}, \ldots, g_{n}$ in 
$ (\O/\h)_{\bar{\p}}$ generate  $(\g/\h)_{\bar{\p}}$ as a complete
intersection ideal of  $ (\O/\h)_{\bar{\p}}$. }
\end{itemize}
Then
$\int _{\h}\g=\g ^*:=(g_1, \ldots , g_t)+(g_{t+1}, \ldots, g_{n})^2  .$
\end{Proposition}

\begin{proof}
 It is obvious that $\g ^*\subset  \int _{\h}\g$. It was proved in
 \cite{JS} that $\O/\g$ and   $\O/\int _{\h}\g$  have the same
 associated ideals. By Lemma
 \ref{property} and  \ref{Ass2}, neither 
 $\g ^* $ nor  $\int _{\h}\g$ has embedded primes
 and these two ideals have common radical $\g$. We only need to prove
they are equal locally at any  associated prime $\p$ of $\O/\g$. 
 
Let $f=a_1g_1+\cdots +a_ng_n \in \int _{\h}\g$. Then for any
 $\xi \in \hder {\h} {\O}$, $\xi(f)\equiv \xi (f')
\equiv 0 \mod  \g$, where $f':=a_{t+1}g_{t+1}+\cdots+a_ng_n $.
By assumption
 $\p \in \cV (\h)\setminus \cV(\J(\h)) $,
$\tilde{\xi} (\tilde{f'})$  is
contained in  $(\g/\h)_{\bar{\p}}$, 
where  we use $\tilde{e}$ to denote the image of an
element $e\in \O$ (resp.  $\hder {\h}{\O}$)
 under taking modulo $\h$ and localization at $\bar{\p}$.
 In other words,
$$\tilde{f'}=\tilde{a}_{t+1}\tilde{g}_{t+1}+\cdots+\tilde{a}_n\tilde{g}_n
\in\left(\left(\int _{\h}{\g}\right)/{\h}\right)_{\bar{\p}}
=\int _{0}\left({\g}/{\h} \right)_{\bar{\p}}
=(\tilde{g}_{t+1},\ldots ,\tilde{g}_n)_{\bar{\p}}^2,$$
where the last equality follows from Lemma~\ref{property} and \cite{Seibt}. 
By the assumption (2), we have 
$\left((a_{t+1},\ldots , a_n)/\h\right)_{\bar{\p}}
\subset \left((g_{t+1},  \ldots ,g_n)/\h\right)_{\bar{\p}},$ 
which implies that
$(a_{t+1},\ldots , a_n)_{\p}\subset \g _{\p}.  $
Thus, we proved that 
$(\g ^* )_{\p}=\left(\int _{\h}{\g}\right)_{\p}$.
\end{proof}

\begin{Corollary}\label{basic2}\sl
Under the  assumption of the Proposition~\ref{basic},
$N$ is a free    $\O/\g$-module with the images 
$\hat{g}_{t+1}, \ldots, \hat{g}_{n}$ of $g_{t+1}, \ldots, g_{n}$ as basis, and
	  $M\simeq T(M)\oplus N$.
\end{Corollary}
\begin{proof}
We only need to prove that $\hat{g}_{t+1}, \ldots, \hat{g}_{n}$
are $\O/\g$-linearly independent. Suppose there are some
$\bar{\beta} _{t+j}\in \O/\g$ such that 
$\hat{a}:= \bar{\beta}_{t+1}\hat{g}_{t+1}+\cdots +\bar{\beta}_{n}\hat{g}_{n}
=0$ in $N$. By the expression of the primitive ideal in
Proposition~\ref{basic} and  taking representatives, we have
$a\in (g_{t+1}, \ldots, g_n)\cap \int _{\h}\g
=(g_1, \ldots , g_t)\cap (g_{t+1}, \ldots, g_n) +(g_{t+1}, \ldots, g_n)^2.$
Let $a\equiv \beta _1g_1+\cdots + \beta _tg_t \mod (g_{t+1}, \ldots, g_n)^2$.
Then $-\beta _1g_1-\cdots - \beta _tg_t+  \beta _{t+1}g_{t+1}+\cdots + 
 \beta _n g_n \equiv 0 \mod \g^2$, from  which and the assumption on
 $\g$, it follows that  $\beta _j \in \g.$   
\end{proof}

\vskip 1mm

Let ${\mathfrak  g}$ be generated by an
 $\O$-regular sequence: 
$g_1, \ldots , g_n. $  Assume that there exist an integer $0\leq t\leq n$ and 
 non-zero divisors
 $\bar{\beta} _1, \ldots ,\bar{ \beta} _t \in \O/\g$ such that
  $\bar{\beta} _1\hat{g}_1, \ldots , \bar{ \beta} _t\hat{g}_t$ are
zero in $M$ as an $\O/\g$-module. Namely
$\beta _1g_1, \ldots , \beta _tg_t\in {\mathfrak  g}^2+{\mathfrak  h}$,
where  $\beta _ig_i$ is the representative of  $\bar{\beta} _i\hat{g}_i$.
 
Let $\{h_1, \ldots,  h_p\}$ be a minimum generating set of 
${\mathfrak h}$. Denote 
$$
h=\transpose (h_1, \cdots , h_p),\quad    g=\transpose(g_1, \cdots ,g_n),\quad
 G=\transpose(G_1, \cdots,  G_t),\quad   \Lambda =\text{diag}\{\beta _1,
 \cdots, \beta _t\},$$
where T means the transposition of the matrix indicated.

Let  $A$ and $B=(B_1\vdots  B_2)$ be the matrices  such that
$\Lambda\transpose( g_1, \cdots,   g_t)=Ah+G ,\,  h=Bg, $
 where $A$ is a $t\times p$ matrix, $B$ a $p\times n$ matrix,  $B_1$ 
a $p\times t$ matrix, $B_2$  a $p\times(n-t)$ matrix, 
 and  $ G_i \in {\mathfrak  g}^2$. 
Let $C_1=AB_1, C_2=AB_2$, then we have
$$(\Lambda-C_1)\transpose(g_1, \cdots , g_t)
-C_2
\transpose(g_{t+1}, \cdots ,  g_n)
\equiv 0 \mod{{\mathfrak  g}^2}$$
 
Note that $\g/\g^2$ is a free $\O/\g$-module.
We have 
$\bar{\Lambda}=\bar{C_1}, \, \bar{C_2}=0 \text{ in } \O/\g.$
From this we obtain the following lemma similar to the implicit function
 theorem.
 
\begin{Lemma}\label{implicitfunctiontheorem}\sl
 Let  $\g $ be a radical  complete intersection ideal of $\O$, then  
$\det \bar{C}_1=\det\bar{\Lambda} =\bar{\beta} _1  \cdots 
\bar{\beta} _t$ is a non-zero divisor in $\O/\g$
and consequently
$\rank (B_1)\ge t$ and $t\le p.$\QED

\end{Lemma}

\section{Freeness of $N$  and the primitive ideal}\label{free}

Let $\h \subset \g $ be complete intersection ideals of
  $\O=k[x_0,x_1,\ldots, x_m]$,
  a polynomial ring over a field $k$ of characteristic zero.
Assume that 
the Jacobian $\J(\h)$ is not contained in any minimal prime of $\O/\g$.
Let $p:=\grade \h, \, n:= \grade \g .$
If   ${\mathfrak g}$ is radical and  generated by an ${\mathcal O}$-regular 
sequence $g_1 , \ldots , g_n$, we can choose the minimal generating set
$\{h_1, \ldots, h_p\}$ of
 ${\mathfrak h}$ such that 
(with changing of the generators of ${\mathfrak g}$ if necessary):
$$h_i\equiv\sum _{j=1}^t b_{ij}g_j  \mod \g ^2,
\quad\quad 1\le i\le p \eqno{(\ref{free}.1)}$$
where   $t$ is an integer
$0\leq t\leq n$,
$b_{ij}\notin {\mathfrak g}\setminus 0$, and   
for each $j$, $(b_{1j}, \ldots , b_{pj})\ne 0$ in $(\O/\g)^p$,
(otherwise one could lower $t$). Denote  $B:=(b_{ij})$.

\begin{Lemma}\label{representation}\sl
Under the assumptions above, we have
\begin{itemize}
\item[(1)]{ $t\geq p$;}
\item[(2)]{There exists at least one maximal minor of $B$ which is 
	  non-zero divisor in $\O/\g$.}
\end{itemize}
\end{Lemma}

\begin{proof}
 Since $\h$ is a 
complete intersection,  the Jacobian $\J(\h)$ of $\h$ can be 
  generated by $\h$ and 
the $p\times p$ minors of the Jacobian matrix $J(\h)$ of
 $\h$. Since each
of these minors, say $\Delta _{j_1,\ldots , j_p}$, 
 is the determinant of $BG_{j_1, \ldots , j_p}$ modulo
${\g}$, where 
$G_{j_1, \ldots , j_p}$ is the $t\times p$ submatrix of   $J(\g)$, 
 consisting of the $0\leq j_1<\cdots j_p\leq m$ columns of 
$J(\g)$, the Jacobian matrix  of $\g$.
 Suppose $t<p$, there would be, 
$\det(BG_{j_1, \ldots , j_p})\equiv 0  \mod  \g$.
This is impossible since we assume that
 $\J(\h)$ is not contained in any minimal prime of $\O/\g$. This proves (1).

(2) Suppose that all the $p\times p$ minors of $B$ are zero divisor in
 $\O/\g$. Then there exists $0\ne a\in\O/\g $ such that
$ab_1\wedge b_2\wedge \cdots \wedge b_p =0 \text{ in } 
 \bigwedge ^p(\O/\g)^t$,
where $b_i\in (\O/\g)^t$ is the image  of the $i$-th row vector of $B$.
 Hence $ab_1, b_2, \ldots , b_p$ are linearly dependent in $(\O/\g)^t$.
 Then there are $a _1, \ldots , a _p\in \O/\g$ which are not all zero, 
such that
$a _1 b_1+\cdots +a _p b_p=0\in (\O/\g)^t.
$
Hence
$a _1 h_1+\cdots +a _p h_p
\equiv 0 \mod  \g ^2 .
$
From this we have
$\J(\h) \subset   \g$, a contradiction.\end{proof}

 \begin{Proposition}\label{free1}\sl
 Let $\h\subset \g$ be   complete intersection ideals of a
  polynomial ring $\O$ over a field $k$ of characteristic zero.
Assume that $\spec{\O/\g}$ is reduced  and connected,
 and  the Jacobian $\J(\h)$ is not contained 
in any minimal prime of $\O/\g$.
  If in {\rm (\ref{free}.1)} we have $t=p=\grade \h$,
 then $b:=\det(b_{ij})$ is  a non-zero divisor in $\O/\g$, and
\begin{itemize}
\item[1)]{ the images $\hat{g}_1, \ldots ,\hat{g}_p$ of
$g_1, \ldots , g_p$ generate $ T(M)$  over $\O/\g$;}
\item[2)]{ the images $\hat{g}_{p+1}, \ldots, \hat{g}_{n}$ of 
$g_{p+1},\ldots, g_n$ generate
 $N$  freely over ${\mathcal O}/{\g}$, so $M\simeq T(M)\oplus N$,
 and $\text{rank} (M)=\text{rank}(N)=\dim X-\dim \Sigma=n-p$;}
\item[3)]{ For each $\p \in \Ass(\O/g)$, the images
$\tilde{g}_{p+1}, \ldots, \tilde{g}_{n}$ of
$g_{p+1}, \ldots , g_{n}$ form an
 $\left(\O/\h\right)_{\bar{\p}}$-regular sequence  and generate
 $\bar{\g}_{\bar{\p}}$;}
\item[4)]{ $ \int _{\h} {\mathfrak  g}
=(g_1, \ldots , g_p)+(g_{p+1},  \ldots , g_n)^2;$}
\item[5)]{ there is a length formula if it is finite
$$
\lambda  (\h,\g) :=l_{\O/\g}(T(M))=l_{\O/\g}\left(
\frac{\mathcal O}{(b)+{\mathfrak g}}\right).
$$}
\end{itemize}
\end{Proposition}
We call $\lambda  (\h, \g)$ the
{\it torsion number } of the pair $(\h, \g)$. When $\h$ and $\g$ are
 clear from the context, we write $\lambda  $ for 
 $\lambda  (\h,\g)$. 

  \begin{proof}
 Since $t=p$ and $b$ is a non-zero divisor,
  one can see that $\hat{g}_1, \ldots ,\hat{g}_p\in  T(M)$ by multiplying
 $B^*$ to the both sides of (\ref{free}.1) , where $B^*$ is the adjoint matrix
 of $B$.

Since $\hat{g}_1, \ldots ,\hat{g}_n$ generate $M$ over $\O/\g$ and 
(\ref{torsionsection}.1)
is exact, $\pi(\hat{g}_{p+1}), \ldots , \pi(\hat{g}_{n})$ generate $N$.
If there is a relation:
$\bar{\beta}_{p+1}\pi(\hat{g}_{p+1})+\cdots+
 \bar{\beta}_n\pi(\hat{g}_{n})=0\in N,$
then
$\bar{\beta}_{p+1}\hat{g}_{p+1}+\cdots+
 \bar{\beta}_n\hat{g}_{n}\in T(M).$
This means that there is a non-zero divisor $\bar{\beta}\in \O/\g$
such that
$\bar{\beta}(\bar{\beta}_{p+1}\hat{g}_{p+1}+\cdots +
\bar{\beta}_n\hat{g}_{n})=0 \in M.$
By taking representatives, this simply means
$\beta\beta _{p+1}g_{p+1}+\cdots +
\beta\beta _n g_{n} \in {\mathfrak g}^2+{\mathfrak h}.$
Hence there are $\mu _1, \ldots \mu _p \in {\mathcal O}$ such that
$$\mu _1 h_1+\cdots + \mu _p h_p+\beta\beta _{p+1}g_{p+1}+\cdots +
\beta\beta _n g_{n} \in {\mathfrak g}^2.$$
By (\ref{free}.1), this becomes
$$\mu '_1g_1+\cdots + \mu '_pg_p+\beta\beta _{p+1}g_{p+1}+\cdots +
\beta\beta _n g_{n} \in {\mathfrak g}^2,$$
where
 $\begin{pmatrix}\mu '_1&\cdots & \mu '_p\end{pmatrix}=
\begin{pmatrix}\mu _1&\cdots  & \mu _p\end{pmatrix}B.$
Since $g_1, \ldots , g_n $ form an ${\mathcal O}$-regular sequence, we have
$\bar{\beta}\bar{\beta }_j=0 $ in $\O/\g$. Note that 
$\bar{\beta}$ is a non-zero divisor, hence $\bar{\beta }_j=0 $ in $\O/\g$.
 This proves 1) and  2).

 For each prime $\p\in \Ass(\O/\g)$,
$N_{\bar{\p}}$
is also free with the images of $\hat{g}_{p+1}, \ldots , \hat{g}_{n}$ as basis.
Then $\left(\int_{0} \bar{\g}\right)_{\bar{\p}} =  \bar{\g}_{\bar{\p}}^2$ since
$\bar{\g} _{\bar{\p}}$ is a reduced  complete intersection in
 $(\O/\h)_{\bar{\p}}$
and  $\p$ is in the regular locus of $\O/\h$ by the assumption.
 Since $(\O/\h)_{\bar{\p}}$ is regular, by Vasconcelos' Theorem \cite{Va1},
$\tilde{g}_{p+1}, \ldots , \tilde{g}_{n}$  is an 
$(\O/\h)_{\bar{\p}}$-regular sequence, and they generated 
$\bar{\g} _{\bar{\p}}$ by Nakayama lemma;

4) follows from  Proposition \ref{basic}.

For the length formula, note that
$$
T(M)=\frac{\int _{\h} {\mathfrak g}}{{\mathfrak g}^2+ {\mathfrak h}}\cong
\frac{(g_1, \ldots, g_p)}{ (g_1, \ldots, g_p)^2+ (g_1, \ldots, g_p) (g_{p+1}, \ldots, g_n)+(h_1,\ldots ,  h_p)}.
$$
It is easy to see that  
$$
M_1:=\frac{(g_1, \ldots, g_p)}{ (g_1, \ldots, g_p)^2+ (g_1, \ldots, g_p) (g_{p+1}, \ldots, g_n)}
$$
is a free $\O/\g$-module. Since $b$ is a non-zero divisor, the following sequence is exact
$$
0\longrightarrow M_1
{\buildrel{\phi _B}\over{\longrightarrow}} M_1
\longrightarrow T(M)\longrightarrow 0,
$$
where $\phi _B(\bar{g}_i):=\sum\limits_{j=1}^p\bar{b}_{ij}\bar{g}_j$. 
By  \cite[A.2.6]{fulton}, we have the length formula of $T(M)$.
\end{proof}

Note that in  the following , we do not assume (\ref{free}.1).

\begin{Proposition}\label{free2}\sl
Let $\h\subset \g$ be complete intersection ideals of $\O$, a 
polynomial ring over a field $k$ of characteristic zero. Let $\grade  \h=p$
and $\grade  \g=n$. 
Assume that $\spec{\O/\g}$ is reduced  and connected,  and $\J(\h)$ 
is  not contained in  any minimal prime of $\O/\g$. 
If $N$ is a free $\O/\g$-module, then 
\begin{itemize}
\item[1)]{there exists an ${\mathcal O}$-regular sequence $g_1, \ldots , g_n$, 
generating ${\mathfrak g}$, such that 
\begin{itemize}
\item{the images $\hat{g}_1, \ldots ,\hat{g}_p$  of
 $g_1, \ldots , g_p$ generate $T(M)$;}
\item{the images $\hat{g}_{p+1}, \ldots ,\hat{g}_n$  of
 $g_{p+1}, \ldots , g_n$ form a basis of $N$;}
\item{rank$(M)$=rank$(N)=n-p=\dim X-\dim \Sigma$;}
\end{itemize}
}
\item[2)]{$\int _{\h}{\mathfrak  g}=(g_1, \ldots , g_p)+(g_{p+1},  \ldots , g_n)^2.$}
\item[3)]{we can choose the generators $h_1, \ldots ,h_p$ of $\h$
such 
{\rm (\ref{free}.1)} holds with  $t=p$ and $b$ a non-zero divisor in $\O/\g$;}
\end{itemize}
\end{Proposition}

\begin{proof}
Let the images $\hat{g}_{t+1}, \ldots ,\hat{g}_n$ of
 $g_{t+1}, \ldots,  g_{n}\in \g$ 
generate $N$ over $\O/\g$, where $t:=n-\text{rank} N$.

For any $\p \in \Ass (\O/\g)$, 
by the assumption, $\p $ is in the regular locus of $\O/\h$, and 
 $N_{\bar{\p}}$ is again a free
module with the images of $\hat{g}_{t+1}, \ldots ,\hat{g}_n$  as basis. 
By   Vasconcelos' theorem (cf. \cite{Va1}), the images
 $\tilde{g}_{t+1}, \ldots, \tilde{g}_n$ of  
$g_{t+1}, \ldots, g_n$ in $(\O/\h)_{\bar{\p}}   $
 form an $(\O/\h)_{\bar{\p}}$-regular sequence. And 
$$\bar{\g}_{\p}=(\tilde{g}_{t+1}, \ldots, \tilde{g}_n)+
 \left(\frac{\int _{\h}{\g}}{\h}\right)_{\bar{\p}}\eqno{(\ref{free}.2)}$$
 Hence the images of $h_1, \ldots, h_p,  g _{t+1}, \ldots, g_n$ in
$\O_{\p} $ form an 
$\O_{\p} $-regular sequence, where  $h_1, \ldots , h_p$ form a minimal
generating set  of $\h$. 
However, since $\h\subset \g \subset \p$, we have
$n-t+p =  \grade (h_1, \ldots, h_p ,  g_{t+1}, \ldots, g_n)_{\p}
\le  \grade (\g)_{\p}=n.$
Hence $t\geq p$.

Extend $g_{t+1}, \ldots , g_n$ to an   ${\mathcal O}$-regular sequence:
 $g_1, \ldots , g_n$, such that they generate ${\mathfrak g}$.
Then  $\hat{g}_1, \ldots , \hat{g}_n$ generate $M$ over $\O/\g$.
 We look for the generating set of $T(M)$. Let
$$
\pi (\hat{g}_i)=\bar{c}_{it+1}\pi(\hat{g}_{t+1})+\cdots +
\bar{c}_{in}\pi(\hat{g}_{n}), \quad  i=1, \ldots , t.$$
Hence by (\ref{torsionsection}.1),
$
\hat{g}'_i:=-\hat{g}_i+ \bar{c}_{it+1}\hat{g}_{t+1}+\cdots +
\bar{c}_{in}\hat{g}_{n}\in T(M), 
 i=1, \ldots , t.
$
Taking representatives, denote
$  g'_i:=-g_i+ c_{it+1}g_{t+1}+\cdots +
c_{in}g_{n},  i=1, \ldots , t, $  $
g'_{t+j}=g_{t+j}, \ j=1, \ldots , n-t.$
Then ${\mathfrak g}=(g'_1, \ldots, g'_n)$, with 
$\hat{g}_1', \ldots, \hat{g}_t'\in T(M)$. By Lemma
\ref{implicitfunctiontheorem},
 $t\leq p$. 
We have proved 1).

Note that $\bar{\g}_{\bar{\p}}$ is also a complete intersection ideal in 
the regular ring $\left(\O/\h\right)_{\bar{\p}}$, and we have $p=t$ in 
(\ref{free}.2). By actually \cite{Seibt} and \cite{hochster},
$\left(\frac{\int _{\h}{\g}}{\h}\right)_{\bar{\p}}=\bar{\g}_{\bar{\p}}^2$.
By Nakayama lemma and (\ref{free}.2), we have  
$\bar{\g}_{\bar{\p}}=(\tilde{g}_{p+1}, \ldots, \tilde{g}_n).$
By Proposition~\ref{basic}, we have 2).

Since ${\mathfrak h} \subset {\mathfrak g}$, we have
$h_i=b_{i1}g'_1+\cdots +b_{ip}g'_p+b_{ip+1}g'_{p+1}+\cdots +b_{in}g'_n,
i=1, \ldots, p.$
For any $\xi \in \hder {\h}{\O}$, we have
$b_{ip+1}\xi(g'_{p+1})+\cdots +b_{in}\xi(g'_n)\equiv 0 \mod{\g}, 
i=1, \ldots, p.$ 
Hence $b_{ip+1}g'_{p+1}+\cdots +b_{in}g'_n\in \int _{\h}{\g}$, which implies
that $b_{ip+1}, \ldots, b_{in} \in {\mathfrak g}$ for $i=1, \ldots, p.$
It is obvious that $b$ is a non-zero divisor in $\O/\g$.
\end{proof}

\begin{Remark}\label{critera}\rm Combining the conclusions in Corollary~\ref{basic2},
Proposition~\ref{free1} and \ref{free2}, one sees that  the Main Theorem 
 is proved.
 Moreover,  either of the  equivalent conditions in the Main Theorem
is equivalent to 3) in Proposition~\ref{free2}.
\end{Remark}

\begin{Example}\label{typical}\rm Let $\h$ be defined by $h:=x^3+xy^3+2x^2z+2z^2=0$, $\g$ be
defined by $g_1:=x^2+y^3=0,g_2:=z=0$. Thus ${\mathfrak h}=(h), {\mathfrak g}
=(g_1, g_2)$.  Notice that 
$h$ is not weighted homogeneous. So it is not easy to find the generator 
set of $\hder {\h}{\O}$. Then we have the same problem for
 $\int_{\h}{\g}$.
If we denote
$g'_1=g_1+2xg_2+g_2^2$, then $h=xg_1'+(2-x)g_2^2$, where $x$ is a non-zero
divisor in $\O/\g$. By Proposition \ref{free1}, we have:
\begin{itemize}
\item{$T(M)$ is generated by $g_1'$ over $\O/\g$}
\item{$\int _{\h} {\mathfrak g}=(g_1', g_2^2)=(x^2+y^3+2xz, z^2)$ }
\item{$N=(g_2)/(g_1'g_2, g_2^2)$ is a free  $\O/\g$-module. }
\end{itemize}
\end{Example}

The following example shows that it is not necessary  for 
$T(M)$ to be generated by $\bar{g_i}$  when $t>p$.

\begin{Example}\label{example0}\rm 
Let  ${\mathfrak g}=(g_1,g_2)$
 with $g_1=xy, g_2=z$ and  ${\mathfrak h}=(h)$ with 
$h=x^2y+yz +z^2=xg_1+yg_2+g_2^2$.
Then $\O/\g \cong {\mathbb C}\{x, y\}/(xy)$, 
$\int _{\h}{\mathfrak g}=(x^2y,yz,z^2)$ (see Example~\ref{example1} for this
 formula) and 
 ${\mathfrak g}^2+{\mathfrak h}=(x^2y^2, xyz, z^2, x^2y+yz)$.
So $T(M)\cong {\mathbb C}x^2y$. And $N$ is not a free $\O/\g$-module.
\end{Example}

\section{Lines on spaces with isolated complete intersection singularities}
\label{application}

We include some applications of the theory to lines on a variety
with isolated complete intersection singularity.
Let $\O$ be the polynomial ring $\cc[x, y_1, \ldots, y_n]$ or the
convergent power series  $\cc\{x, y_1, \ldots, y_n\}$.
Let $\Sigma$ be a line in ${\mathbb C}^{n+1}$ defined by $\g$.
Define $_{\Sigma}{\mathcal K}:={\mathcal R}_{\Sigma}\rtimes {\mathcal C}$,
 the semi-product of ${\mathcal R}_{\Sigma}$ with  the contact group 
${\mathcal C}$ (cf. \cite{mather}), where ${\mathcal R}_{\Sigma}$ is a
 subgroup of ${\mathcal R}:={\rm Aut}(O)$ consisting of all the 
$\varphi \in {\mathcal R}$  preserving $\g$.
This group has an action on the space
 ${\mathfrak m}{\mathfrak g}{\mathcal O}^p.$  For
 $h=(h_1, \ldots, h_p)\in {\mathfrak m}{\mathfrak g}{\mathcal O}^p$,
there is an ideal $\h$ generated by
$h_1, \ldots, h_p$, and a variety 
$X={\mathcal V}({\mathfrak h})$.  The image of the differential of $h$:
${\mathcal O}^{n+1}{\buildrel {dh^*}\over {\longrightarrow}}
{\mathcal O}^{p}$
is denoted by ${\rm th}(h)$. 
Define a $_{\Sigma}{\mathcal K}$-invariant
$$\tilde{\lambda}:=\tilde{\lambda} (\Sigma , X)
=\dim_{\mathbb C}\frac{{\mathcal O}^p}
{{\rm th}(h)+{\mathfrak g}{\mathcal O}^p  }.
$$
Choose $\Sigma $ as the $x$-axis.
Then $\Sigma$ can be defined by ${\mathfrak g}=(y_1, \ldots, y_n)$.
Denote $\O _X : \O/\h, \O _{\Sigma}:=\O/\g.$

\begin{Proposition}\label{linesonicis}\sl Let $\Sigma$ be a line on 
a variety  $X$ with isolated complete intersection singularity of
 codimension $p$. Then $h$ is 
$_{\Sigma}{\mathcal K}$-equivalent to an $\tilde{h}$ with  components
$\tilde{h}_i\equiv b_{i}y_i \mod \g ^2$,
where
$b_{i}\notin {\mathfrak   g}$, $i=1, \ldots p$. Moreover
$\lambda (\g, \h)=\tilde{\lambda} (\Sigma, X)=\dim_{\mathbb C}
\O /
(b+\g) =\sum_{k=1}^pl _i.$
where $b:=b_{1}\cdots b_{p}$, and $l_i$ is the valuation of $\bar{b}_{i}$
in $\O _{\Sigma}$.
\end{Proposition}
\begin{proof} 
Since
 $\Sigma \subset X$,
  for a given generating set $\{h_1, \ldots, h_p\}$ of ${\mathfrak h}$, 
we have
$h_i\equiv \sum\bar{b}_{ij}y_j \mod{{\mathfrak   g}^2 },  i=1,
 \ldots , p,$
where $ \bar{b}_{ij}\in \O _{\Sigma}$, and for fixed $i$, $ \bar{b}_{ij}$'s
are not all zero since  $X$ is complete intersection and
 $X_{\text{sing}}=\{0\}\subsetneq \Sigma $.
 Since $ \O _{\Sigma}$ is a principal ideal domain, by changing the indices,
 we can assume that $\bar{b}_{11}\mid \bar{b}_{ij}$. Let
$y'_1=y_1+\sum_{j=2}^n \frac{\bar{b}_{1j}}{\bar{b}_{11}}y_j.$
Then
$h_1\equiv \bar{b}_{11}y_1' \mod  {\mathfrak   g}^2 .$
Let 
$h_i'=h_i-\frac{\bar{b}_{i1}}{\bar{b}_{11}}h_1,  i=2, \ldots , p.$
Repeat the above argument will prove the first part of the proposition.

Consider the exact sequence
$$
{\mathcal O}^{n+1} \quad {\buildrel {dh^*}\over {\longrightarrow}} \quad
{\mathcal O}^{p}\longrightarrow \text{ coker}(dh^*)\longrightarrow 0.
$$
By tensoring with $\O _{\Sigma}$, we have the  exact sequence
$$
\O _{\Sigma}^{n+1}\quad
 {\buildrel {\bar{dh^*}}\over {\longrightarrow}}\quad
\O _{\Sigma}^{p}\longrightarrow\text{ coker}(\bar{dh^*})\longrightarrow 0.
$$
However by the expression of $h_i$'s above, this is just
$$
\O _{\Sigma}^{p}{\buildrel {\bar{dh^*}}\over {\longrightarrow}}
\O _{\Sigma}^{p}\longrightarrow
{\frac{{\mathcal O}^p}{{\rm th}(h)+{\mathfrak g}{\mathcal O}^p  }} \longrightarrow 0.
$$
Since $ \bar{b}\ne 0$, by \cite[A.2.6]{fulton}, we have the formula for
$\tilde{\lambda}$.
\end{proof}

\begin{Corollary}\label{linecase}\sl
Let  $X$ be a variety  with isolated complete intersection singularity of
 codimension $p$  in ${\mathbb C}^{n+1}$,  and
 $\Sigma $  a line in $X$ defined by   ${\mathfrak g}$.
 Then we can choose the coordinates of  ${\mathbb C}^{n+1}$ such that
 ${\mathfrak g}  =(y_1, \ldots, y_n)$, 
 $\hat{y}_1, \ldots , \hat{y}_p\in T(M)$ and 
$\hat{y}_{p+1}, \ldots , \hat{y}_n$ generate $N$ which is free of rank
 $n-p$, and 
$$\int _{\h} {\mathfrak  g}=(y_1, \ldots , y_p)+(y_{p+1},  \ldots , y_n)^2.
\eqno{\square}
$$
\end{Corollary}

\begin{Example}\label{example1}\rm  Consider the situation in 
Example \ref{example0}. Denote
${\mathfrak g}_1=(x, z)$ and  ${\mathfrak g}_2=(y,z)$. They define 
$\Sigma _1=y$-axis
and $\Sigma _2=x$-axis respectively. We have 
${\mathfrak g}={\mathfrak g}_1\cap {\mathfrak g}_2$
 and $\Sigma =\Sigma _1\cup \Sigma _2\subset X$.
Since $h=(y+z)z+yx^2$, by Corollary \ref{linecase},
 $\int _{\h} {\mathfrak g}_1=(x^2, z)$.
Since $h=(x^2+z)y+z^2$, again  by Corollary \ref{linecase},
 $\int _{\h}{\mathfrak g}_2=(y, z^2)$.
These tell us that
 $\int _{\h}{\mathfrak g}
=\int _{\h}{\mathfrak g}_1\cap \int _{\h}{\mathfrak g}_2$.
\end{Example}

\begin{Corollary}\label{orollary2}{\rm  (Due to Pellikaan) } \sl
Let  $X$ be a variety  with isolated complete intersection singularity of
 codimension $p$  in ${\mathbb C}^{n+1}$, and
 $\Sigma $  a line in $X$, defined by   ${\mathfrak g}=(y_1, \ldots, y_n)$.
Then  the Second Exact Sequence \cite{matsumura, hartshorne} is exact on 
 the left also:
$$
0\longrightarrow M
 \longrightarrow
\Omega  ^1_X\otimes {\mathcal  O}_\Sigma\longrightarrow 
\Omega ^1_\Sigma \longrightarrow 0.
$$
 Furthermore it is splitting and 
$$
T(\Omega  ^1_X\otimes {\mathcal  O}_\Sigma)=T(M),\quad
\rank(\Omega  ^1_X\otimes {\mathcal  O}_\Sigma)=n-p-1.
$$
These tell us that the torsion number $\lambda (\h,\g)$ is independent of the 
choice of the generating sets of ${\mathfrak g}$ and ${\mathfrak h}$.
\end{Corollary}

 \begin{proof} We have the following 
presentation:
$$
{\mathcal O}_X^{p}{\buildrel{dh}\over{\longrightarrow}}
{\mathcal O}_X^{n+1} \longrightarrow \Omega  ^1_X \longrightarrow 0.
$$
Tensoring with $\O _{\Sigma}$, we have the exact sequence
$$
\O _{\Sigma}^{p}{\buildrel{\overline{dh}}\over{\longrightarrow}}
\O _{\Sigma}^{n+1} \longrightarrow \Omega  ^1_X\otimes {\mathcal  O}_\Sigma
 \longrightarrow 0.
$$
Remark that the map $\overline{dh}$ is equivalent to a map defined by the 
matrix $(\bar{b}_{ij})$. Hence 
$$
\Omega  ^1_X\otimes {\mathcal  O}_\Sigma
\cong \frac{\O _{\Sigma}^{n+1}}{\text{im}\overline{dh}} 
\cong \O _{\Sigma}^{n-p+1}\oplus   
{\frac{\O _{\Sigma}^{p}}{\text{im}(\bar{b}_{ij})}}
\cong \O _{\Sigma}\oplus N\oplus T(M)\cong \O _{\Sigma}\oplus M.
$$
Then exact is
$$
0\longrightarrow M\longrightarrow 
{\frac{\O _{\Sigma}^{n+1}}{\text{im}\overline{dh}} }
\longrightarrow \O _{\Sigma}\longrightarrow 0.
$$
Since $\Omega ^1_\Sigma$ is free  $\O _{\Sigma}$-module of rank 1, by 
 \cite[A.2.2]{fulton},  we have the  exact sequence.
\end{proof}

\begin{Remark}\label{lambda}\rm
In general, given an analytic space germ $(X, 0)$, one cannot find a
 smooth curve $L\subset X$ that passes through and is  not contained in
 $X_{\rm sing}$. However, if there are smooth curves on $X$ in the above 
 sense, how to distinguish them is a problem. We found that the torsion
number $\lambda$ is a nice candidate for this purpose \cite{JS1}.
In  studying the Euler-Poincar\'{e} characteristic $\chi(F)$ of the
 Milnor fibre $F $ of 
 a function with singular locus a smooth curve on a singular space,
 we found that this $\lambda$ also appears in $\chi(F)$.
Note also that the torsion number was generalized to 
`` higher torsion numbers'' in \cite{JS, JS1}.

\end{Remark}

\vskip 1mm

\noindent{\bf Acknowledgments } This article comes from part of 
the author's
thesis which  was finished at Utrecht University under the advising of
 Professor Dirk Siersma. Many discussions with R. Pellikaan, A. Simis,
J.  R. Strooker
and H. Vosegaard were very helpful. The author thanks all of them.


\begin{thebibliography}{99}
\bibitem{AM}
{\sc M. F. Atiyah, I. G. MacDonald, } 
{\it Introduction to commutative algebra},    Addison-Wesly Publ. Comp. 1969.


\bibitem{fulton}
{\sc  W. Fulton, } {\it Intersection theory,}
 Ergeb. Math. Grenzgebiete, 3. Folge 2, Springer-Verlag, Berlin, 1988.

\bibitem{hartshorne}
{\sc  R. Hartshorne, } {\it  Algebraic geometry,} Graduate texts 
in mathematics, {\bf 52}, Springer-Verlag, 1977.


\bibitem{hochster}
{\sc M. Hochster, }  Criteria for equality of ordinary and symbolic powers of
 primes, {\it Math. Z.}, {\bf 133} (1973) 53-65.

\bibitem{J}
{\sc G. Jiang, }  Functions with non-isolated singularities on singular spaces, 
Thesis, Universiteit Utrecht, 1997.

\bibitem{JS}
{\sc G. Jiang, A.   Simis, }  Higher relative primitive ideals, 
Tokyo Metropolitan University mathematics preprint series 1999 no 10,
to appear in: {\it Proc. Amer. Math. Soc.}


\bibitem{JS1}
{\sc G.  Jiang,  D.  Siersma, }  Local embeddings of lines in singular
 hypersurfaces,  {\it  Ann. Inst. Fourier (Grenoble)},
 {\bf 49} (1999) no. 4, 1129--1147.


\bibitem{mather}
{\sc  J. N. Mather, }  Stability of $C^{\infty}$-mappings III: Finitely
determined map germs, {\it Inst. Hautes \'Etudes Sci. Publ. Math.,} \ 
 {\bf 35}(1968)  127-156.


\bibitem{matsumura}
{\sc  H. Matsumura, } {\it Commutative algebra,} The Benjamin/Cummings,
 Reading, 1980.

\bibitem{P1}
{\sc R. Pellikaan, } Hypersurface singularities and resolutions of Jacobi
 modules, Thesis, Rijkuniversiteit Utrecht, 1985.

\bibitem{P2}
{\sc  R.  Pellikaan, } Finite determinacy of functions with non-isolated
 singularities, {\it  Proc. London Math. Soc.}, (3){\bf 57} (1988) 357--382.

\bibitem{Saito}
{\sc  K. Saito, } Theory of logarithmic differential forms and logarithmic 
vector fields, {\it J. Fac. Sci. Univ. Tokyo Sect. 1A Math.}
 {\bf 27} (1980) 265--291.

\bibitem{Seibt} 
{\sc P. Seibt, } Differential filtrations and symbolic powers
of regular primes, {\it Math. Z.} {\bf 166} (1979) 159--164.

\bibitem{Aron1}
{\sc A. Simis, }
 Effective computation of symbolic powers by jacobian matrices,
 {\it Comm. in Algebra,} {\bf 24} (1996) 3561--3565.

\bibitem{Va1}
 {\sc  W. V. Vasconcelos, } Ideals generated by $R$-sequences, {\it J. Alg.,}
{\bf 6} (1967) 309-316.

\bibitem{vas2} 
{\sc W. V. Vasconcelos, } {\it Computational Methods in
Commutative Algebra and Algebraic Geometry}, Algorithms and Computation
in Mathematics, Vol. {\bf 2}, Springer-Verlag, 1998.
\end{thebibliography}
\end{document}